\documentclass[12pt]{article}
\usepackage{amsmath,amssymb,theorem}
\newtheorem{thm}{Theorem}%%[section]

\newtheorem{cor}{Corollary}%[section]
\newtheorem{lem}{Lemma}%[section]
\theorembodyfont{\rmfamily}
\newtheorem{rem}{Remark}%[section]
\newenvironment{proof}{\noindent{\bf Proof.}}{\hfill
  $\blacksquare$\par\noindent}

\def\bbe{{\mathbb E}}

\def\bbp{{\mathbb P}}
\def\bbr{{\mathbb R}}
\begin{document}
\title{Concentration for Infinitely Divisible Vectors 
with Independent Components}

\author{C. Houdr\'e\thanks{Laboratoire d'Analyse et de Math\'ematiques
Appliqu\'ees, CNRS UMR 8050, Universit\'e Paris XII, 94010 Cr\'eteil Cedex,
France and School of Mathematics, Georgia Institute of Technology,
Atlanta, GA 30332, USA.}
\and 
P. Reynaud-Bouret\thanks{DMA Ecole Normale Sup\'erieure, 45 rue  
d'Ulm 75230 Paris Cedex 05, France and School of Mathematics, Georgia
Institute of Technology, Atlanta, GA 30332, USA.}}

\maketitle
\begin{abstract}
For various classes of Lipschitz functions we provide dimension free 
concentration inequalities for 
infinitely divisible random vectors with independent components and 
finite exponential moments.
\end{abstract}

\baselineskip=20pt
The purpose of this note is to further visit the concentration phenomenon
for infinitely divisible vectors with independent 
components in an attempt to obtain dimension free
concentration.

Let $X\sim ID(\gamma,0,\nu)$ be an infinitely divisible (i.d.) vector (without Gaussian 
component) in
$\bbr^d$, and with characteristic function 
$\varphi (t)=\bbe e^{i\langle t,X\rangle}$, $t\in\bbr^d$ (throughout,
$\langle\cdot,\cdot\rangle$ denotes the Euclidean inner product
in $\bbr^d$, while $\|\cdot\|$ is the corresponding Euclidean norm).  As well
known,
\begin{equation}\label{eq1}
\varphi (t)=\exp\left\{i\langle t,\gamma\rangle +\int_{\bbr^d} (e^{i\langle t,u\rangle}
-1-i\langle t,u\rangle{\bf 1}_{\|u\|\le 1})\nu (du)\right\}, 
\end{equation}
where $\gamma\in\bbr^d$ and where $\nu\not\equiv 0$ (the L\'evy measure)
is a positive Borel measure on $\bbr^d$, without atom at the origin
and such that $\int_{\bbr^d} (1\wedge \|u\|^2)\nu(du)<+\infty$.  As also well known, 
$X$ has independent components if and only if
$\nu$ is supported on the axes of $\bbr^d$, i.e.,
\begin{equation}\label{defindpdc}
\nu(dx_1,\dots, dx_d)=\sum^d_{k=1} \delta_0 (dx_1)\cdots \delta_0(dx_{k-1})\tilde\nu_k
(dx_k)\delta_0(dx_{k+1})\cdots \delta_0 (dx_d).
\end{equation}
Moreover, the independent components of $X$ have same law if 
and only if, the one dimensional L\'evy measures $\tilde\nu_k$ are
the same measure denoted by $\tilde\nu$. 

Below, and throughout, by  $f$ Lipschitz with 
constant $a$ we mean that $|f(x) - f(y)|\le a\|x-y\|$, 
for all $x, y \in \bbr^d$
(the Lipschitz convention stated in \cite{h} also applies).   
Let us start by recalling the following simple lemma which will be crucial to our
approach \cite{hpas}. 

\begin{lem}\label{simple}
Let $X\sim ID(\gamma,0,\nu)$ be such that $\bbe \|X\|^2<+\infty$.
Let $f,g:\bbr^d\to \bbr$ be Lipschitz functions.  Then,
\begin{align}\label{eq2}
\bbe f(X)g(X)&-\bbe f(X)\bbe g(X)\nonumber\\
&=\int^1_0 \bbe _z\left[ \int_{\bbr^d}(f(U+u)-f(U))(g(V+u)-g(V))\nu(du)\right]dz,
\end{align}
where the expectation $\bbe_z$ is with respect to the i.d.\ vector, $(U,V)$
in $\bbr^{2d}$ of parameter $(\gamma,\gamma)$ and with 
L\'evy measure $z\nu_1+(1-z)\nu_0$,
$0\le z\le 1$.  The measure $\nu_0$ is given by
$$\nu_0 (du,dv)=\nu(du)\delta_0(dv)+\delta_0 (du)\nu(dv), u, v \in \bbr^d, $$
while $\nu_1$ is the measure $\nu$ supported on the main diagonal of
$\bbr^d$.
\end{lem}

An important feature of the representation \eqref{eq2} is the fact that the
first marginal of $(U,V)$ is $X$ and so is its second marginal.

With the above framework and denoting by $e_1,e_2,\dots, e_d$, the canonical
basis of $\bbr^d$, we first prove:

\begin{thm}\label{thm1}
Let $X=(X_1,\dots, X_d)\sim ID(\gamma,0,\nu)$ have independent components 
and be such that $\bbe e^{t\|X\|}<+\infty$, for some $t>0$.  Let 
$f:\bbr^d\to\bbr$, and let there exist $b_k \in \bbr$, $k=1,\dots, d$, such that
$|f(x+ue_k)-f(x)|\le b_k|u|$, for all $u \in\bbr$, $x\in \bbr^d$.
Let $$h_f(t)=\sup_{x\in \bbr^d}\sum^d_{k=1}\int_{\bbr}|f(x+ue_k)-f(x)|^2\,
\frac{e^{tb_k|u|}-1}
{b_k|u|}\tilde\nu_k(du) \mbox{, } 0\le t<M,$$ where 
$M=\sup\left\{t>0:\forall\ k=1,\dots, d, \bbe e^{tb_k|X_k|}<+\infty\right\}$.  Then
\begin{equation}\label{eq3}
\bbp(f(X)-\bbe f(X)\ge x)\le e^{-\int^x_0 h^{-1}_f(s)ds},\end{equation}
for all $0<x<h_f^{-1}(M^-)$.
\end{thm}

\begin{proof} The proof is akin to proofs given in \cite{h}, and 
the above result complements the results there.
First, by independence,
\begin{align*}
C&=
\left\{t>0:\forall\ k=1,\dots, d, \bbe e^{tb_k|X_k|}<+\infty\right\}\\
&= \left\{t>0:\forall\ k=1,\dots, d,\int_{|u|>1} e^{tb_k|u|}\tilde\nu_k(du) 
< +\infty\right\}.\end{align*}
Next, we apply the covariance representation \eqref{eq2} to $f$ 
satisfying the above hypotheses and moreover assumed to be 
bounded and such that $\bbe f=0$.  Thus,
\begin{align*}
\bbe fe^{tf} &=\int^1_0 \bbe _z\left[e^{tf(V)}\sum^d_{k=1}
\int_{\bbr}(f(U+ue_k)-f(U))
(e^{t(f(V+ue_k)-f(V))}-1)\tilde\nu_k(du)\right]dz\\
&\le\int^1_0 \!\bbe_z\!\!\left[\!e^{tf(V)}\sum^d_{k=1}\int_{\bbr}|f(U+ue_k)-f(U)|
|f(V+ue_k)-f(V)|\,\frac{e^{tb_k|u|}-1}{b_k|u|} \,\tilde\nu_k(du)\right]dz\\
&\le\!\!\int^1_0\!\bbe_z\left[e^{tf(V)}\!\sum^d_{k=1}\!\int_{\bbr}\!\frac{
|f(U+ue_k)-f(U)|^2\!
+\!|f(V+ue_k)-f(V)|^2}{2}
\!\!\left(\!\!\frac{e^{tb_k|u|}-1}{b_k|u|}\!\!\right)\!\,\tilde\nu_k(du)\!\right]\!dz\\
&\le h_f(t)\bbe\left[ e^{tf}\right],
\end{align*}
where we have used the ``marginal property" mentioned above and since
$h_f(t)$ is well defined for $0\le t<M$.  Integrating this last inequality, 
applied to $f-\bbe f$, leads to
\begin{equation}\label{eq4}
\bbe e^{t(f-\bbe f)}\le e^{\int^t_0 h_f(s)ds},\qquad 0\le t<M,
\end{equation}
for all $f$ bounded satisfying the hypotheses of the theorem.  
Fatou's lemma allows to remove the boundedness assumption in
\eqref{eq4}.

To obtain the tail inequality \eqref{eq3}, the Bienaym\'e-Chebyshev
inequality gives
$$
\bbp (f(X)-\bbe f(X)\ge x)\le \exp\left(-\sup_{0<t<M}\left(tx-\int^t_0
h_f(s)ds\right)\right)=e^{-\int^x_0 h^{-1}_f(s)ds},$$
by standard arguments, e.g., see \cite{h}.
\end{proof}

Theorem~\ref{thm1} is a bit formal, and we are now going to provide various
cases where more concrete estimates are possible.  Our first
corollary, of Bennett--Prokhorov type, 
improves the constants in a result of \cite{h}.  If the components of $X$ are 
iid Poisson random variables,
then \eqref{eq5} recovers also a result obtained by Bobkov and
Ledoux \cite{bl2} via modified log-Sobolev inequalities.  This corollary is 
optimal in the one dimensional case, but suboptimal in the multidimensional one 
(see Corollary~\ref{cor6}, for a more dimension free result)  

\begin{cor}\label{cor1}
Assume the hypotheses of the previous theorem.  Moreover, let 
$\nu$ have bounded support with
$$R_k=\inf\{\rho>0: \tilde\nu_k(|x|>\rho)=0\}.$$
Set $bR=\max_{1\le k\le d} b_kR_k$ and set
${\bar a}^2=\sup_{x\in\bbr^d}\sum^d_{k=1}\int_{|u|\le R_k}|f(x+ue_k)-f(x)|^2\tilde\nu_k(du)$.
Then, for all $x>0$,
\begin{equation}\label{eq5}
\bbp (f(X)-\bbe f(X)\ge x)\le 
e^{-\frac{{\bar a}^2}{b^2R^2}\ell\left(\frac{bRx}{{\bar a}^2}\right)},
\end{equation}
where $\ell (u)=(1+u)\log (1+u)-u$, $u>0$.  
\end{cor}

\begin{proof} It is enough to note that $M=+\infty$ and that
$$h_f(s)\le {\bar a}^2\left(\frac{e^{sbR}-1}{bR}\right).$$
Integrating the reciprocal, gives $\frac {-x}{bR}+
\left(\frac x{bR}+\frac{{\bar a}^2}{b^2R^2}\right)\log\left(1+\frac{bRx}{{\bar a}^2}
\right)$.  
\end{proof}

Let us now give a result which holds for L\'evy
measures with unbounded support, 
giving a Bernstein type inequality.

\begin{cor} \label{cor2}
Assume the hypotheses of Theorem~\ref{thm1}. Let $X\sim
ID(\gamma,0,\nu)$ have iid components and let  $f$ be  such that  $\tilde a^2=
\sup\limits_{\buildrel{\scriptstyle x\in\bbr^d}\over
      {\scriptstyle u\in\bbr, u\ne 0}}
\sum^d_{k=1}\frac{|f(x+ue_k)-f(x)|^2}{|u|^2}$, and 
$b=\max_{1\le k\le d} b_k$ are finite.   Then for all
$0<xb/\tilde a^2<h^{-1}(M^-)$, 
\begin{equation}\label{eq9}
\bbp (f(X)-\bbe f(X)\ge x)\le 
\exp^{\left(-\frac{\tilde a^2}{b^2}\int^{xb/\tilde a^2}_0 h^{-1}(s)ds\right)},
\end{equation}
where $h(s)=\int_{\bbr} |u|(e^{s|u|}-1)\tilde\nu(du)$. 

Moreover, if there exist
$C>0$ and $V^2>0$ such that
\begin{equation}\label{eq6}
\int_{\bbr} |u|^n\tilde\nu(du)\le \frac{C^{n-2}n!}2 \, V^2,\qquad
\forall n\ge 2,
\end{equation}
then, for all $x>0$, 
\begin{equation}\label{eq7}
\bbp(f(X)-\bbe f(X)\ge x)\le e^{-\frac{\tilde a^2V^2}{b^2C^2}\,\ell\left(
\frac{bCx}{\tilde a^2V^2}\right)},\end{equation}
where now $\ell (u)=(1+u)-\sqrt{1+2u}, u>0$. 
\end{cor}

\begin{proof} Again, we just need to bound $h_f$ of Theorem~\ref{thm1}.
For \eqref{eq9}, we bound $h_f$ by 
$$h_f(s) \le \frac{\tilde{a}^2}{b}\int_\bbr u^2
(e^{sbu}-1)\tilde{\nu}(du),$$
and the result follows.

The condition \eqref{eq6} implies exponential moments, for 
$0<t<\frac1{bC}$.  Moreover, 
\begin{align*}
h_f(t)&\le \sup_x\int_{\bbr}\sum^d_{k=1}\frac{|f(x+ue_k)-f(x)|^2}{|u|^2}
\frac{|u|(e^{tb|u|}-1)}b\,\tilde\nu(du)\\
&\le \frac{\tilde a^2}b \frac{V^2}{2C}\sum^\infty_{k=2} k(tbC)^{k-1}\\
&=\frac{\tilde a^2}b\frac{V^2}{2C}\left(\frac1{(1-tbC)^2}-1\right),
\end{align*}
using \eqref{eq6}.  Integrating its reciprocal, we get
$$\bbp(f(X)-\bbe f(X)\ge x)\le e^{-\frac{\tilde a^2V^2}{b^2C^2}
\left(1+\frac{bCu}{\tilde a^2 V^2}-\sqrt{1+\frac{2bCu}{\tilde a^2V^2}}
\right)}.$$
\end{proof}

\begin{rem}\label{rem0}
(i) An instance of the potential suboptimality of the
previous results is the case of the (symmetric) exponential measure.
Indeed, if $X_1,\dots, X_d$ are iid with density $2^{-1}{e^{-|x|}}$,
then the exponent in \eqref{eq7} or in \eqref{eq9} is of order
$\min\left(\frac xb\,,\frac{x^2}{\tilde a^2}\right)$, while an
inequality of Talagrand \cite{t} asserts that the order
$\min\left(\frac xb\,,\frac{x^2}{a^2}\right)$, where
\begin{equation}\label{defa}
a^2=\sup_{x,u\in \bbr^d}\frac{|f(x+u)-f(x)|^2}{\|u\|^2},
\end{equation} 
holds true.  Clearly, $a^2\le \tilde a^2\le db^2$.  It is then clear
that \eqref{eq7} or \eqref{eq9} are optimal for linear functions, or
infimum like above but not for the Euclidean norms.  Actually, the
example of the Euclidean norm, i.e., $f(x)=\|x\|$, for which $a^2=1$,
$b^2=1$, while $\tilde a^2=db^2=d$ shows that the concentration
inequalities obtained that way are not dimension free. This is after
all quite natural since Theorem~\ref{thm1} is really a result about
``$\ell_1$-Lipschitz'' functions.  We could also (to mimic the L\'evy
measure of the exponential law) replace in \eqref{eq6} $n!$ by
$(n-1)!$, but the corresponding estimate will not be dimension free
either.   Our next result will show that for the Euclidean norm, a better estimation
of \eqref{eq3} leads to dimension free concentration (see also
Theorem~\ref{thm2} and Theorem~\ref{thm3}).  

(ii) Here is, however, an example of a class of function for which we
can exactly get a dimension-free exponential inequality.
Let $X$ be an iid vector as in Corollary \ref{cor2}, for which 
$\tilde{\nu}$ has a support included  in $\bbr^+$. Let $f$ be defined by
$$f(x)= \inf_{\alpha \in A} f_\alpha(x),$$
where the $\{f_\alpha,\alpha\in A\}$ are non decreasing coordinates by
coordinates and such that for all $1\le k\le d$ and $\alpha\in A$, there
exists a constant $b_{\alpha,k}$ such that
$$\forall u \in \bbr^+, \forall x \in \bbr^d, 0 \le f_\alpha(x+ue_k)-
f_\alpha(x) \le b_{\alpha,k} u.$$
Then we can apply Corollary \ref{cor2} to $f$. First let us  remark that the
supremum in $\tilde{a}$ can be taken for $u>0$, and that
$$\tilde{a}^2\le \sup_{x\in \bbr^d, u>0} \sum_{k=1}^d
\frac{|f_{\hat{\alpha}}(x+ue_k)-f_{\hat{\alpha}}(x)|^2}{u^2},$$
where $\hat{\alpha}$ is the index where $\inf_{\alpha\in A}f_\alpha(x)$ is achieved.
Therefore, $\tilde{a}^2\le \sup_{\alpha\in A} \sum_{k=1}^d
b_{\alpha,k}^2$.  On the other hand, one also has $b= \sup_{1\le k\le d,
\alpha\in A} b_{\alpha,k}$.  This gives $f$, Lipschitz, nonlinear, 
with $\tilde{a}< \sqrt{d} b$.  In particular, $f(X) = \inf_{1 \le k
\le d} X_k$ verifies these conditions with $\tilde{a}=a=b=1$, and Corollary
\ref{cor2} implies a complete dimension-free  exponential inequality
for the deviations of $f$ above its mean.

\end{rem}

\begin{cor}\label{cor3}
Let $X\sim ID(\gamma,0,\nu)$ have independent components and be such that
$\bbe e^{t\|X\|}<+\infty$, for some $t>0$.
Let $M=\sup\{t>0:\forall\ k=1,\dots, k$, $\bbe e^{t|X_k|}<+\infty\}$.
Let $\varepsilon >0$.
Then, for all $0<x<h(M^-)$
\begin{equation}\label{eq10}
\bbp (\|X\|\ge (1+\varepsilon)\bbe \|X\|+x)\le e^{-\int^x_0
h^{-1}(s)ds},\end{equation}
and
\begin{equation}\label{eq10bis}
\bbp (\|X\|\le (1-\varepsilon)\bbe \|X\|-x)\le e^{-\int^x_0
h^{-1}(s)ds},\end{equation}
where the (dimension free) function $h$ is
given by $h(t)=8\max_{1\le k\le d}\int_{\bbr}|u|(e^{t|u|}-1)\tilde\nu_k(du)
+\frac2{(\varepsilon\bbe \|X\|)^2}\sum^d_{k=1}\int_{\bbr}|u|^3
(e^{t|u|}-1)\tilde\nu_k(du)$.
\end{cor}

\begin{proof} We apply Theorem~\ref{thm1} to $f(x)=(\|x\|-\varepsilon\bbe \|X\|)^+$.
Let us compute the various parameters and integrals for this $f$.
First, it is easily verified that for each $k$, 
$|f(x+ue_k)-f(x)|\le\left|\|x+ue_k\|-\|x\|\right|
{\bf 1}_{A_k}$, where  
$A_k=\{\|x+ue_k\|\ge \varepsilon \bbe \|X\|$ or $\|x\|\ge \varepsilon \bbe \|X\|\}$, 
and where $\varepsilon > 0$.   
We then have 
\begin{equation}\label{eq11}
|f(x+ue_k)-f(x)|\le\frac{|u(2x_k+u)|{\bf 1}_{A_k}}{\|x+ue_k\|+\|x\|}
\le \frac{2|ux_k|}{\|x\|}+\frac{u^2}{\varepsilon\bbe \|X\|},   
\end{equation}
where $x_k$ is the kth coordinate of $x$.  
Moreover, since $|f(x+ue_k)-f(x)|\le |u|$, we have
\begin{align*}
&\sum^d_{k=1}\int_{\bbr}|f(x+ue_k)-f(x)|^2\frac{e^{tb_k|u|}-1}{b_k|u|}\,
\tilde\nu_k(du)\\
&\qquad \le \sum^d_{k=1}\int_{\bbr}\left( 8u^2\frac{|x_k|^2}{\|x\|^2}
+\frac{2u^4}{(\varepsilon\bbe \|X\|)^2}\right)\left(\frac{e^{t|u|}-1}{|u|}\right)
\tilde\nu_k(du).\end{align*}
Hence $h_f$ in Theorem~\ref{thm1} is such that
$$h_f(t)\le 8\max_{1\le k\le d}\int_{\bbr}|u|(e^{t|u|}-1)\tilde\nu_k(du)+
\frac2{(\varepsilon\bbe \|X\|)^2}\sum^d_{k=1}\int_{\bbr}|u|^3
(e^{t|u|}-1)\tilde\nu_k(du).$$
To finish the proof of (\ref{eq10}) note that
$\|X\|-\varepsilon\bbe \|X\|\le (\|X\|-\varepsilon\bbe \|X\|)^+$ and that 
$\bbe (\|X\|-\varepsilon\bbe \|X\|)^+\le \bbe \|X\|$.  
To get the lower bound (\ref{eq10bis}), just proceed as
above but with the function $f(x)=-(\|x\|-\varepsilon\bbe \|X\|)^+$
and note that $(\|X\|-\varepsilon\bbe \|X\|)^+\le \|X\|$ and that 
$(1-\varepsilon)\bbe \|X\| \le \bbe (\|X\|-\varepsilon\bbe \|X\|)^+$.

\end{proof}

\begin{rem}\label{rem1}
(i) The function $h$ in the previous result is dimension free.
Indeed,
$$h(t)\le 8\max_{1\le k\le d}\int_{\bbr}|u|(e^{t|u|}-1)\tilde\nu_k(du)
+\frac{2d}{(\varepsilon\bbe \|X\|)^2}\max_{1\le k\le d}\int_{\bbr}
|u|^3(e^{t|u|}-1)\tilde\nu_k(du),$$
but,
$$d\min_{1\le k\le d}(\bbe |X_k|)^2\le (\bbe \|X\|)^2\le 
d\max_{1\le k\le d}\bbe(X_k^2).$$

(ii) When the L\'evy measure has bounded support, the previous result leads 
under the assumptions of Corollary~\ref{cor1} to

\begin{equation}\label{eq11a}
\bbp (\|X\|\ge (1+\varepsilon)\bbe 
\|X\|+x)\le e^{-\left(\frac{x}{R}+\frac{V^2}{R^2}\right) 
\log\left(1+\frac{xR}{V^2}\right) + \frac{x}{R}},
\end{equation} 
where 

\begin{equation}\label{eq11ab}
V^2= 8 \max_{1\le k \le d}\left( \int_{|u|\le R} u^2 \tilde\nu_k(du) \right)+ 
\frac{2}{(\varepsilon\bbe \|X\|)^2}
\sum_{k=1}^d \int_{|u|\le R} u^4 \tilde\nu_k(du), 
\end{equation}
does not depend on the dimension $d$.  This implies that there exists
a constant $C>0$ independent of $d$ such that
\begin{equation}\label{eq11aa}
\bbp (\|X\|\ge (1+\varepsilon)\bbe \|X\|+x) \le e^{- C
\min\left(\frac{x}{R}\log\left(\frac{xR}{V^2}\right), \frac{x^2}{V^2}\right)},
\end{equation} 
for all
$x>0$. (\eqref{eq11a} and \eqref{eq11ab} allow to improve, for the Euclidean
norm, the range
and the constants in the last theorem of \cite{hm}.)

A direct consequence of Corollary 1 in \cite{h} is the fact that for
$X$ infinitely divisible in $\bbr^d$ with boundedly supported L\'evy
measure,
\begin{equation}\label{eq11b}
\bbe e^{\frac{\|X\|}{R} \log^+\left(\frac{\lambda \|X\|}{R}\right)} <
+\infty, \end{equation} 
for all $\lambda >0 $ such that $\lambda V^2/R^2 < 1/e$, where $V^2 =
\int_{\|u\|\le R}\|u\|^2\nu(du)$.  Although tight (take a one
dimensional Poisson random variable with mean one), \eqref{eq11b} is
not optimal.  Indeed a result of Rosi\'nski \cite{r} asserts that
(for i.d.~vectors in Banach spaces) \eqref{eq11b} holds under the
tighter condition $\lambda p_0 < 1/e$, where $p_0 = \nu(\|u\|=R)$.
Similarly, another direct consequence of \eqref{eq11a} is the
following fact.
\begin{cor}\label{rosinski}
Let  $X$ be as in Corollary~\ref{cor3} above,
\begin{equation}\label{eq11bc}
\bbe e^{\frac{\|X\|}{R} \log^+\left(\frac{\lambda \|X\|}{ R}\right)} <
+\infty,
\end{equation}
for all $\lambda >0 $ such that $\lambda V^2/R^2 < 1/e$, where now
$V^2$ is given by \eqref{eq11ab}.
\end{cor}  
Hence, for vectors with independent
components, this last condition on $\lambda$ is dimension free and in
this sense, it improves on the general result obtained in \cite{r}.
Although dimension free, the condition $\lambda V^2/R^2 < 1/e$ (with
$V^2$ as in \eqref{eq11ab}) is not optimal (again, take a one
dimensional Poisson random variable with mean one).  In view of
\cite{r}, and say for $X$ with iid components, one might wonder if 
$\lambda p_0 < 1/e$, where $p_0 = {\tilde\nu}(|u|=R)$ might be optimal.  

The estimate (\ref{eq11aa}) also improves a case of the exponential
inequality derived for suprema of integrals with respect to a centered
inhomogeneous Poisson process in \cite{rb}. Let $N$ be a Poisson
process on $\mathbb{X}$ with intensity $s$ with respect to $\mu$. Let
$\cal P$ be a partition of $\mathbb{X}$ and $S$ the space of piecewise
constant functions on $\cal P$. Let 
$$\chi = \sup_{f\in S}\frac{\int_{\mathbb{X}} f \frac{dN-s
d\mu}{\mu(\mathbb{X})}} {\sqrt{\int_{\mathbb{X}} f^2 \frac{d\mu}{\mu(\mathbb{X})}}}.$$
For this special choice of $S$, Proposition 9 of \cite{rb} implies 
that, for all positive $\varepsilon$,
there exists $C>0$, such that
\begin{equation}\label{eq11bcd}
\forall x>0, \bbp(\chi\ge (1+\varepsilon) \sqrt{\bbe \chi^2} + x)
\le \exp\left[-C \min\left((\sqrt{\eta \mu(\mathbb{X}) })x,
\frac{\mu(\mathbb{X})}{K}x^2\right)\right], 
\end{equation}
where $\eta = \inf_{I\in \cal P} \mu(I)$ and $K=\sup_{I\in \cal P}
\frac{\int_I s d\mu}{\mu(I)}.$

\noindent But $\chi$ can be viewed as the Euclidean norm of the infinitely
divisible vector
$$X=\left(\int_I\frac{dN-sd\mu}{\sqrt{\mu(I)\mu(\mathbb{X})}}\right)_{I\in
\cal P},$$ with independent components. The L\'evy measures of the
components (see (\ref{defindpdc})) are given by
$$\tilde{\nu}_I = \left(\int_I s d\mu\right)
\delta_{\frac{1}{\mu(I)\mu(\mathbb{X})}}.$$  
Thus, we can apply  
(\ref{eq11a}) or (\ref{eq11aa}) with $R=1/\sqrt{\eta \mu(\mathbb{X}) }$ and
$$V^2= c(\varepsilon)\frac{K}{\mu(\mathbb{X})}.$$
Above, the constant $c(\varepsilon)$ does not depend on $\cal P$ or $\mu({\mathbb{X}})$ as soon
as $\eta>1$ which is the interesting case where this type of
inequality leads to adaptive estimator of the intensity $s$.  We refer
to \cite{rb} for a complete description of this procedure.
Therefore, (\ref{eq11aa}) gives an extra logarithmic factor with
respect to (\ref{eq11bcd}) when $S$ is a space
of piecewise constant functions on a given partition.  More precisely, 
for all positive $\varepsilon$,
there exists $C>0$, such that
\begin{equation}\label{eq11bce}
\forall x>0, \bbp(\chi\ge (1+\varepsilon) {\bbe \chi} + x)
\le \exp\left[-C \min\left(\sqrt{\eta \mu(\mathbb{X}) }
x\log\left(\frac{x\sqrt{\eta }\mu(\mathbb{X})^{3/2}}{K}\right),
\frac{\mu(\mathbb{X})}{K}x^2\right)\right]. 
\end{equation}

(iii) If $X$ has iid components $X_1,\dots, X_d$ and if $X_1$ has an
exponential distribution with density $2^{-1}e^{-|x|}$, $x \in \bbr$,
and L\'evy measure $|u|^{-1}e^{-|u|}$, $u\in \bbr, u\neq 0,$ the
previous result is a version of Talagrand's inequality (for norms).
Indeed, in this case $M=1$ and since $\bbe(|X_k|)=1$, we obtain for
all $0\le t<1$,
\begin{eqnarray*}
h(t) &=& 8 \int_{\bbr} |u| (e^{t|u|}-1) \frac{e^{-|u|}}{|u|} du 
+ \frac{2}{\varepsilon^2 (\bbe 
\|X\|)^2} d \int_{\bbr} 
|u|^3 (e^{t|u|}-1) \frac{e^{-|u|}}{|u|} du, \\
&\le & 16 \int_0^\infty (e^{tu}-1) e^{-u} du 
+ \frac{4}{\varepsilon^2} \int_0^\infty u^2 (e^{tu}-1) e^{-u} du, \\
& \le & 16 \left(\frac1{1-t} -1\right) + \frac{8}{\varepsilon^2} 
\left(\frac1{(1-t)^3} -1\right),\\
& \le & \left(16 + \frac{8}{\varepsilon^2}\right) \left(\frac{1}{(1-t)^3}-1\right).\\
\end{eqnarray*}
 
This leads to 
$$\forall x>0,~~\bbp (\|X\|\ge (1+\varepsilon)\bbe \|X\|+x) \le e^{- x -
\frac{3}{2}\left(16 + \frac{8}{\varepsilon^2}\right) \left(1- \left( 1 +
\frac { x} {\left(16 + \frac{8}{\varepsilon^2}\right)}\right)^{2/3}\right)},$$
which implies that
$$\forall x>0,~~\bbp (\|X\|\ge (1+\varepsilon)\bbe \|X\|+x)
\le  \exp\left(-\frac{x^2}{6(16+\frac{8}{\varepsilon^2})+4 x}\right).$$  
If one is only interested  in the order of magnitude of the deviation
of $\|X\|$, this is completely equivalent to Talagrand's inequality 
applied to the Euclidean norm, since 
(forgetting the constants and the dependency in $\varepsilon$) 
the exponent above is of order $-\min(x,x^2)$.  
However, one may want to get the exact upper deviation
of $\|X\|$ from its mean (and not a constant times its mean).  
To see the difference, let us look at the reverse form :
$$\forall u>0,~~\bbp (\|X\|\ge (1+\varepsilon)\bbe
\|X\|+\frac{4\sqrt{3u}}{\varepsilon}+4\sqrt{6u}+4u)\leq e^{-u}.$$
We can then minimize in $\varepsilon$ and get:
$$\forall u>0,~~\bbp (\|X\|\ge \bbe
\|X\|+4(3u)^{1/4}(\bbe\|X\|)^{1/2}+4\sqrt{6u}+4u)\leq e^{-u}.$$
But $\bbe \|X\|$  grows like $\sqrt{d}$.  So for $d$ large, the quadratic term
disappears and this is equivalent to
$$\forall x>0,~~\bbp (\|X\|\ge \bbe
\|X\|+x)\leq e^{-C\min(x, \frac{x^4}{d})},$$
for some constant $C$. 
Hence, for $d$ large, our method loses the 
quadratic behavior with respect
to \cite{t}.

Our result is more restrictive than Talagrand's since it is only
proved for norms rather than for arbitrary Lipschitz functions, and
cannot give the exact order for the upper deviations from the mean but it
is also more general since valid for any i.d.\ law with finite
exponential moments (note too that the Lipschitz image of the exponential can be a bounded 
random variable and thus not i.d. and that not any i.d. variable with
exponential moment is a
Lipschitz image of the exponential variable. 

(iv) A generalization of Corollary \ref{cor3} to
$\|X\|_A=\sqrt{X^*AX}$, where $A = (a_{j,k})$ is a symmetric positive
definite matrix, is also possible.  It is sufficient to remark that for all $x$ in $\bbr^d$, 
$\|Ax\|^2 \le \lambda_{\mbox{\tiny max}} \|x\|_A^2$ where
$\lambda_{\mbox{\tiny max}}$ is the largest
eigenvalue of $A$. Then we can apply Theorem \ref{thm1} to 
$f(x) = (\|x\|_A - \varepsilon \bbe \|X\|_A)^+$, noticing that
$b_k=\lambda_{\mbox{\tiny max}}$ works and that
$$h_f(t) \le 8\lambda_{{\mbox{\tiny max}}}\max_{1\le k\le d}
\int_{\bbr}|u|(e^{t|u|}-1)\tilde\nu_k(du)
+\frac2{(\varepsilon\bbe \|X\|_A)^2}\sum^d_{k=1}a_{k,k}^2\int_{\bbr}|u|^3
(e^{t|u|}-1)\tilde\nu_k(du).$$
This upper bound is dimension free since $ \|x\|_A^2 \ge
\lambda_{\mbox{\tiny min}}
\|x\|^2$ where $\lambda_{\mbox{\tiny min}}$ is the smallest eigenvalue of $A$.
\end{rem}

We can in fact prove a result true for every Lipschitz function, 
by using the same type of method.

\begin{thm}\label{thm2}
Let $X$ be as in {\rm Theorem~\ref{thm1}}.  Let $f:\bbr^d\to\bbr$ be
Lipschitz, with constant $a$.  Then, 
$$\bbp \left(f(X)\ge \bbe f(X)+ a\sqrt{2\sum^d_{k=1}
{\rm Var}\, X_k} + ax \right)\le e^{-\int^x_0 h^{-1}(s)ds},
$$
for all $0<x<h(M^-)$, where now 
$$h(t)=8\max_{1\le k\le d}\int_{\bbr}|u|(e^{t|u|}-1)\tilde\nu_k(du)
+\frac2{\sum^d_{k=1} {\rm Var}\, X_k}\sum^d_{k=1}\int_{\bbr}|u|^3
(e^{t|u|}-1)\tilde\nu_k(du).$$
\end{thm}

\begin{proof} We apply Theorem~\ref{thm1} to $\phi (X)=\sqrt{\bbe_Y\|X-Y\|^2},$
where $Y$ is a vector such that $\bbe_Y\|Y\|^2 <+\infty$ and independent 
of $X$.    
As, $\sqrt{\bbe_Y\|\cdot\|^2}$ is a norm (for vectors depending
on $Y$), we have
$$|\phi (X+ue_k)-\phi (X)|\le\sqrt{
\bbe_Y\|X+ue_k-Y-(X-Y)\|^2}\le |u|,$$
Thus $b_k=1$, for all $k =1,\dots, d$.  Also, 
\begin{equation}\label{eq15}
|\phi (X+ue_k)-\phi (X)|^2=
\left( \frac{\bbe_Y(2u(X_k-Y_k)+u^2)}{
\sqrt{\bbe_Y\|X-Y\|^2} + \sqrt{\bbe_Y\|X+ue_k-Y\|^2}}\right)^2.
\end{equation}
Note that $\phi (X)\ge\sqrt{\sum^d_{k=1}{\rm Var}\, Y_k}$.  Hence,
the right hand side of \eqref{eq15} is dominated by 
\begin{equation}\label{eq16}
\frac{8u^2\bbe_Y (X_k-Y_k)^2}{\bbe_Y\|X-Y\|^2}+\frac{2u^4}{\sum^d_{k=1}{\rm Var}\, Y_k}\,.
\end{equation}
We then see 
(using \eqref{eq16}) that the function $h_\phi$ 
in Theorem~\ref{thm1} is such that
$$h_\phi(t)\le 8\max_{1\le k\le n}\int_{\bbr}|u|(e^{t|u|}-1)\tilde\nu_k(du)+
\frac{2\sum^d_{k=1}\int_{\bbr}|u|^3
(e^{t|u|}-1)\tilde\nu_k(du)}{\sum^d_{k=1}{\rm Var}\,Y_k}.$$  

\noindent Returning to $f$, and taking for 
$Y$ an independent copy of $X$, we get 
$$\bbe \phi (X)\le\sqrt{\bbe_X \bbe_Y\|X-Y\|^2}\!=\! \sqrt{2\sum^d_{k=1}{\rm Var}\, X_k}.$$ 
Moreover,    
$|f(X)-\bbe f(X)|\le a\phi (X)$.  These last two estimates finally give 
$$\bbp \left(f(X)\ge \bbe f(X)+ a\sqrt{2\sum^d_{k=1}{\rm Var}\, X_k}
+ ax\right)\le \exp\left(-\int^x_0 h^{-1}(t)dt\right).$$
\end{proof}

\begin{rem}\label{remthm2}
The above result gives a dimension-free exponential rate of decay 
for the deviations of $f(X)$ above its mean plus $a\sqrt{d}$, up to
some constants. For the exponential distribution, Theorem~\ref{thm2} does 
not give an exponential rate with two speeds, one
using $b$ (defined in Corollary \ref{cor2}) and the other using $a$
(defined by (\ref{defa})).  This cannot be seen either in Corollary
\ref{cor3}, since for the Euclidean norm $a=b$.  But one can combine 
Theorem~\ref{thm1} and Theorem~\ref{thm2} together.  For iid variables with
exponential symmetric distributions, this gives a partial version of
Talagrand's result \cite{t}.  First, we look at the deviation of $f$ above
$m=\bbe f(X) +2a\sqrt{d}$.  As $\bbe f(X) \le m$, from Corollary~\ref{cor2}
there exists some absolute constant $c_1$ such that
$$\bbp \left(f(X)\ge m +x\right)\le \exp\left(-c_1\min\left(\frac{x}{b}, \frac{x^2}
{\tilde{a}^2}\right)\right),$$
and from Theorem~\ref{thm2} there exists some absolute constant
$c_2$ such that 
$$\bbp \left(f(X)\ge m +x\right)\le \exp\left(-c_2\min\left(\frac{ x}{a},
\frac{x^2}{ a^2}\right)\right).$$  This implies that there exists some
absolute constant $c_3$ such that
$$\bbp \left(f(X)\ge m +x\right)\le \exp(-c_3 g(x))$$
where 
\begin{eqnarray*}
g(x) & = &\frac{x^2}{a^2}, \mbox{ for } 0\le x\le a,\\
& = &\frac{x}{a}, \mbox{ for } a\le x\le \frac{\tilde{a}^2}{a},\\
& = &\frac{x^2}{\tilde{a}^2}, \mbox{ for } \frac{\tilde{a}^2}{a}\le
x\le \frac{\tilde{a}^2}{b},\\
& = & \frac{x}{b}, \mbox{ for } \frac{\tilde{a}^2}{b}\le x.
\end{eqnarray*}
Thus we recover Talagrand's result for small and large
$x$. In the middle, we have intermediate rate. If $a=b=1$ (as for the Euclidean
norm) or if $a=\tilde{a}$ (as for linear functionals), we  recover exactly 
Talagrand's rate on the whole real line.  
For the deviation with respect to $\bbe f(X)$ and not $\bbe f(X) + 2a\sqrt d$, 
the previous rates become worse, but sometimes improve the
rate given by Corollary~\ref{cor2} for some special parts of
the real line.
\end{rem}

The next result is an easy consequence of Theorem~\ref{thm2} by
applying the same methods as in the proof Corollary~\ref{cor1}.
Combined with Corollary~\ref{cor1}, it will give dimension free rates
in $e^{-x^2/a^2}$, for $x$ small above $a\sqrt d$, and of order
$e^{-\frac{x}{bR}\log x}$, for $x$ large.

\begin{cor}\label{cor6}
Let $X$ be as in Theorem \ref{thm1}.  Moreover, let 
$\nu$ have bounded support with
$$R=\max_{1\leq k\leq d}\inf\{\rho>0: \tilde\nu_k(|x|>\rho)=0\}.$$
Let $f$ be a Lipschitz function with constant $a$.

Then, for all $x>0$,
\begin{equation}\label{lipborne}
\bbp \left(f(X)\ge \bbe f(X)+ a\sqrt{2\sum^d_{k=1}
{\rm Var}\, X_k}+ ax\right)\le 
e^{-\frac{v^2}{R^2}\ell\left(\frac{Rx}{v^2}\right)},
\end{equation}
where $\ell (u)=(1+u)\log (1+u)-u$, $u>0$ and 
$$v^2=8\max_{1\le k\le d}\int_{\bbr}|u|^2\tilde\nu_k(du)
+\frac2{\sum^d_{k=1} {\rm Var}\, X_k}\sum^d_{k=1}\int_{\bbr}|u|^4
\tilde\nu_k(du).$$

\end{cor}
 
\begin{rem}\label{remPoincare}
(i) The above improves Corollary~\ref{cor1} as one can see on a vector of iid
Poisson variables with parameter 1. The quantity $\bar{a}$ appearing in
Corollary~\ref{cor1} is then equal to $\tilde{a}$ appearing in 
Corollary~\ref{cor2} and is of order $\sqrt{d}$, while 
Corollary~\ref{cor6} gives a dimension-free exponential rate of decay
for the deviations of $f$  above   $\bbe f(X)+ a \sqrt{2d}$.

(ii) A natural question is then to know whether or
not the above result is a consequence of, or implies, a result 
of Bobkov and Ledoux
\cite{bl1} which asserts that a Poincar\'e inequality does imply
Talagrand's.  This is not the  case. First, a uniform random variable 
on $[0,1]$ satisfies a
Poincar\'e inequality but is not infinitely divisible.  Second, a Poisson
random variable  has finite exponential moments, is infinitely
divisible but does not satisfy a Poincar\'e inequality. However,
Corollary~\ref{cor6} combined with Corollary~\ref{cor1} gives
dimension free rates in $e^{-x^2/a^2}$, for $x$ small above $a\sqrt d+\bbe f(X)$ and of order
$e^{-\frac{x}{bR}\log x}$, for $x$ large.  This is almost a dimension free inequality with 
two rates except that $f$ has to exceed $\bbe f(X) +a\sqrt{d}$ and not just $\bbe f(X)$ and 
that there are smaller rates for intermediate $x$.

\end{rem}

Of course, we would like a result using  only $a$ and $b$ for every 
Lipschitz functions to exactly recover the
exponential case. In particular, even if $f$ has to exceed a multiple
of $a\sqrt{d}$, we would like to improve the
rates obtained in Remark~\ref{remthm2} when $a\le x \le \tilde{a}^2/b$ and $a > b$.  
The next two results give some further knowledge in this direction.  The first
one deals with concave functions and so also leads to a left tail 
inequality for the Euclidean norm.

\begin{cor}\label{cor4}
Let $X$ be as in {\rm Theorem~\ref{thm1}}, let $f:\bbr^d\to\bbr$ be
concave and let $\tilde b_k=\left|\bbe\, \frac{\partial f(X)}{\partial x_k}
\right|$, $k=1,\dots, d$.  Let 
$M=\sup\left\{t>0:\forall\ k=1,\dots, d, \bbe e^{t\tilde b_k|X_k|}<+\infty\right\}$.  
Let ${\rm Cov} (X, \nabla f(X)) =  \bbe\langle
X,\nabla f(X)\rangle  -  \langle
\bbe X, \bbe \nabla f(X)\rangle $.
Then, 
\begin{equation}\label{eq12}
\bbp (f(X) - \bbe f(X)\ge  -{\rm Cov}(X, \nabla f(X)) +x)\le e^{-\int^x_0 h^{-1}(s)ds},
\end{equation}
for all $0<x<h(M^-)$, where $h$ is given by
$h(t)=\sum^d_{k=1}\int_{\bbr}\tilde b_k|u|(e^{t\tilde b_k|u|}-1)
\tilde\nu_k(du)$, $0 < t < M$.  
\end{cor}

\begin{proof} Since $f$ is concave, and if $Y$ is an independent copy of $X$,
$$f(X)-\bbe f(X)\le \bbe_Y(\langle X-Y,\nabla f(Y)\rangle):=\phi (X).$$
We then apply Theorem~\ref{thm1} to $\phi$.  Indeed,
$$\phi (X+ue_k)-\phi (X)=u\bbe \frac{\partial f(Y)}{\partial x_k}.$$
Hence, $h_\phi (t)=\sum^d_{k=1}\int_{\bbr}\tilde b_k|u|
(e^{t\tilde b_k|u|}-1)\tilde\nu_k(du)$, and the result follows.
\end{proof}

\begin{rem}\label{rem2}
Above, if $|f(x)-f(y)|^2\le a^2\|x-y\|^2$ we get:
\begin{align*}
h_\phi (t)&\le \sum^d_{k=1}\tilde b^2_k\int_{\bbr}|u|^2
\frac{(e^{t\tilde b_k|u|}-1)}{\tilde b_k|u|}\,\tilde\nu_k(du)\\
&\le\max_{1\le k\le n}\int_{\bbr} u^2\frac{(e^{t\tilde b_k|u|}-1)}{\tilde b_k|u|}
\, \tilde\nu_k(du)\sum^d_{k=1}\left(\bbe \frac{\partial f(X)}{\partial x_k}\right)^2
\\
&\le a^2\max_{1\le k\le d}\int_{\bbr}|u|\left(\frac{e^{t\tilde b_k|u|}-1}{\tilde b_k}\right)
\tilde\nu_k(du).\end{align*}
Moreover,
$$0\le \bbe \phi(X)= -{\rm Cov}(X, \nabla f(X))\le  
\bbe\left( \|X-\bbe X\| \|\nabla f(X)\|\right)
\le a \bbe\|X-\bbe X\|.$$
Combining these two facts we see that \eqref{eq12} becomes
\begin{equation}\label{eq13}
\bbp \left(f(X) - 
\bbe f(X) \ge a \bbe \|X-\bbe X\| + x\right)
\le e^{-\int^x_0 h^{-1}(t)dt},\end{equation}
where now $h(t)={a^2}\max_{1\le k\le 1}\int_{\bbr}|u|
\frac{(e^{\tilde b_k|u|}-1)}{\tilde b_k}\tilde\nu_k(du)$.  This last 
inequality is once again dimension free.  
In the particular case, $f(X)=-\|X\|$ ($a=b_k=1$), 
we get $ -{\rm Cov}(X, \nabla f(X)) = \bbe\|X\| - \langle
\bbe X, \bbe \left(X/{\|X\|}\right) \rangle$ and \eqref{eq12} 
becomes      

\begin{equation}\label{eq13bis}
\bbp\left(-\|X\|\ge - \langle
\bbe X, \bbe \left(X/{\|X\|}\right) \rangle  +x\right)\le e^{-\int^x_0 h^{-1}(t)dt},
\end{equation}
with $h(t)=\max_{k=1, \dots, d}\int_{\bbr}|u|(e^{t|u|}-1)\tilde\nu_k(du)$.  
However, the inequality \eqref{eq13bis} does not present any interest when $\bbe X = 0$.
\end{rem}

The second result deals with general Lipschitz functions, gives
exponential inequalities using $a$ and $b$ and allows us to improve
the rates, in the exponential case, when $a\sqrt{d} < \tilde{a}^2/b$
and $a>b$.

\begin{thm}\label{thm3}
Let $X$ be as in Theorem~\ref{thm1}. Let $f$ be a Lipschitz function
with constant $a$, and let $b_k \in \bbr$, $k=1,\dots, d$, such that
$|f(x+ue_k)-f(x)|\le b_k|u|$, for all $u \in\bbr$, $x\in \bbr^d$.
Let $\varepsilon>0$.
Then, for all $0<x<h^{-1}(M)$
\begin{equation}\label{toutesctes}
\bbp\left(f(X)\ge f(0)+a
\bbe\left([\|X\|-\varepsilon]_+\right)+a\varepsilon +x\right) \le 
e^{-\int_0^x h^{-1}(s) ds}, 
\end{equation}
where $$h(t)= 2 a 
\sqrt{\sum_{k=1}^d \left(\int_\bbr |u| (e^{tb_k
|u|}-1)\tilde\nu_k(du)\right)^2 } + \frac{a}{\varepsilon}\sum_{k=1}^d \int_\bbr  |u|^2
(e^{t b_k|u|}-1) \tilde{\nu}_k(du),$$ for all $0\le t < M$.
\end{thm}

\begin{proof}
Let $M>t\ge 0$. First, we have  
$$\bbe\left(\left[f(X)-f(0)- a\varepsilon- a
\bbe\left([\|X\|-\varepsilon]_+\right)\right] e^{tf(X)}\right)\le
\mbox{Cov}(g(X), e^{tf(X)}),$$
where $g(X)= a [\|X\|-\varepsilon]_+$. 
By using equation~(\ref{eq2}), we get that
\begin{multline*}
 \mbox{Cov}(g(X), e^{tf(X)}) =\int^1_0\!\!\bbe_z\!\!\left[e^{tf(V)}\sum^d_{k=1}
\int_{\bbr}\!(g(U+ue_k)-g(U))
(e^{t(f(V+ue_k)-f(V))}-1)\tilde\nu_k(du)\!\right]\!\!dz\\
\le\int^1_0 \bbe_z\left[e^{tf(V)}\sum^d_{k=1}\int_{\bbr}|g(U+ue_k)-g(U)|
|f(V+ue_k)-f(V)|\,\frac{e^{tb_k|u|}-1}{b_k|u|} \,\tilde\nu_k(du)\right]dz.
\end{multline*}
By using the computations done in the proof of Corollary \ref{cor3},
we know that
$$|g(U+ue_k)-g(U)|\le 2a|u| \frac{|U_k|}{\|U\|}+\frac{a
u^2}{\varepsilon}.$$
Let us define $A$ and $B$ by
$$A=  \int^1_0 \bbe_z\left[e^{tf(V)}\sum^d_{k=1}\int_{\bbr} 2a \frac{|U_k|}{\|U\|}
|f(V+ue_k)-f(V)|\,\frac{e^{tb_k|u|}-1}{b_k} \,\tilde\nu_k(du)\right]dz,$$
$$B = \bbe\left[e^{tf(V)}\sum^d_{k=1}\int_{\bbr}\frac{a
|u|}{\varepsilon}
|f(V+ue_k)-f(V)|\,\frac{e^{tb_k|u|}-1}{b_k}
\,\tilde\nu_k(du)\right].$$
Then we obtain that 
$$\bbe\left(\left[f(X)-f(0)- a\varepsilon- a
\bbe\left([\|X\|-\varepsilon]_+\right)\right] e^{tf(X)}\right) \le
A+B.$$
We can bound $A$ by
\begin{align*}
A & \le  \int^1_0 \bbe_z\left[e^{tf(V)}\sum^d_{k=1}\int_{\bbr} 2a \frac{|U_k|}{\|U\|}
|u|\, (e^{tb_k|u|}-1) \,\tilde\nu_k(du)\right]dz \\
&\le 2a
\sqrt{\sum_{k=1}^d \left(\int_\bbr |u| (e^{tb_k
|u|}-1)\tilde\nu_k(du)\right)^2 } \bbe(e^{tf(X)}).
\end{align*}
Similarly, we get the following upper bound for $B$:
$$
B  \le \frac{a}{\varepsilon}\left(\sum_{k=1}^d \int_\bbr |u|^2 (e^{tb_k
|u|}-1)\tilde\nu_k(du)\right)\bbe(e^{tf(X)}).
$$
It remains to use the classic integration/maximisation method  to
conclude the proof.
\end{proof}
\begin{rem}\label{rems}
(i)  Comparing Theorem~\ref{thm2} and Theorem~\ref{thm3}, we see
that $f(0)$ is replacing $\bbe f(X)$. This is not a problem since
(see \cite[Appendix V]{ms}) one can pass from the former to the later up to
some multiplicative constant. In fact, it was already possible to derive
directly Theorem~\ref{thm2} with $f(0)$ instead of $\bbe f(X)$.

(ii) If $X$ is a vector of iid  variables with
density $2^{-1}e^{-|x|}$, and if $b=\max_k b_k$, after
computations (similar to the ones given in Remark~\ref{rem1} (iii)),
we obtain that for every Lipschitz function $f$,  for all $x,\varepsilon>0$,
\begin{equation}
\bbp(f(X)\ge f(0)+ a\sqrt{d}+ a \varepsilon + x ) \le e^{-c\min(\frac{x}{b},
\frac{x^2}{ab\sqrt{d}+2ab \frac{d}{\varepsilon}})},
\end{equation}
for some absolute constant $c$.
The reverse form of this last inequality is more practical in order 
to better understand 
the various orders of
magnitude: for all $x,\varepsilon>0$,
\begin{equation}
\bbp(f(X)\ge f(0)+ a\sqrt{d}+ a \varepsilon +
\square\sqrt{\frac{abdx}{\varepsilon}}+\square \sqrt{ab\sqrt{d}
x}+\square bx ) \le e^{-x},
\end{equation}
where the $\square$ are known absolute constants. So first if one take
$\varepsilon =\delta \sqrt{d}$ then this implies  
$$\bbp\left(f(X)\ge f(0)+ (1+2\delta) a\sqrt{d} + 
\left(\square+\frac{\square}{\delta^2}\right)b x\right) \le e^{-x}.$$
Thus, once $f$ has exceeded $f(0)$ plus a multiple, as close to $1$ as
we want, of $a\sqrt{d}$ the behavior is linear, and the slope is $b$ up to
some multiplicative constant, increasing as $\delta$ tends to $0$. 
This improves the  results of Remark~\ref{remthm2} for the
exponential case when $\tilde{a}^2/b >  > a \sqrt{d}$ and $\tilde{a}>a>b$
since now the linear rate $x/b$ is true on a larger interval.  
This partially recovers Corollary~\ref{cor3} since
$\sqrt{x} \leq \sqrt{d}+x $, for all $d\ge 1$ and $x>0$.

One can also optimize in $\varepsilon$, getting that for all positive $x$
$$\bbp ( f(X) \ge f(0) + a\sqrt{d} + \square a^{2/3} b^{1/3} d^{1/3}
x^{1/3}+ \square a^{1/2} b^{1/2} d^{1/4} x^{1/2}+ \square b x) \le e^{-x}.$$

(iii) This result also improves the rates for iid Poisson variables
with parameter $1$. When Corollary \ref{cor6} gives the rate $\exp(-
C_1 \frac x a \log(\frac x a))$ for the deviations above $Ef+
a\sqrt{d}+x$ for sufficiently large $x$, Theorem \ref{thm3} gives $
  \exp(-
C_1 \frac x b \log(\frac x {a\sqrt(d)}))$ for $x > C_3 a \sqrt{d}$
which is better than Corollary \ref{cor1} and Corollary \ref{cor6} as soon as
$b << a$ and $a\sqrt{d} <<\tilde{a}^2/b$.

(iv) More generally, if one is interested in Lipschitz function of
i.d. vectors with independent components and L\'evy measure with bounded
support, the equivalent of Corollary \ref{cor6} can be obtained by
applying Theorem \ref{thm3} to L\'evy measures with bounded support.
Similarly, the equivalent of Corollary \ref{rosinski} for
$f$,  Lipschitz function with constant $a$ ($b$ being  defined as usual) can also be obtained.   
One straightforward application is then to
say that $|f(X)-f(0)|\le a \|X\|$, giving :
$$\bbe e^{\frac{f(X)}{aR}\log^+(\frac{\lambda f(X)}{aR})} < \infty,$$
for all $\lambda >0$ such that $\lambda V^2/R^2 < 1/e$, where $V^2$ is
given by~(\ref{eq11ab}).
But one may wonder if the above remains true with $b$ instead of $a$, i.e., the
$\ell^1$-Lipschitz constant.  By applying Theorem \ref{thm3} with
$\varepsilon= \sqrt{d}$, it follows that
\begin{equation}
\label{rosinskib}
\bbe e^{\frac{f(X)}{bR}\log^+\left(\frac{\lambda f(X)}{bR}\right)} < \infty,
\end{equation}
for all $\lambda >0$ such that $\lambda \frac{aV^2}{bR} <
1/e$, where this time $$V^2 =3 \sqrt{d} \left(\max_{1\le k \le
d}\left(\int_\bbr |u| \tilde{\nu}_k(du)\right)\vee\max_{1\le k \le
d}\left(\int_\bbr |u|^2 \tilde{\nu}_k(du)\right)\right).$$
As $V^2$ is not dimension free, this is not as sharp as Corollary
\ref{rosinski} for the Euclidean norm, but it is sharper than the 
results of \cite{r} since, in that case, $V^2$ would be of order
$d$. It also implies with $b$ instead of $a$ the following result
$$\bbe e^{\frac{f(X)}{A}\log^+(\frac{ f(X)}{bR})} < \infty,$$
for all $A>bR$, which is a complete dimension free result and which
can be of interest if $a>>b$.
\end{rem}

\

The various results presented here for vectors with finite exponential moment
as well as the general methodology presented in
\cite{hm} delineate quite well the
concentration phenomenon for infinitely divisible vectors.   Nevertheless, and say,
for iid components, it will be interesting to prove versions of Theorem~\ref{thm2}
or of Theorem~\ref{thm3} for the deviations of an arbitrary Lipschitz function  
above its mean and not just above
its mean plus $a\sqrt{d}$, up to a constant.  Such a possible extension would then give,
when combined with Theorem~\ref{thm1} a dimension-free
exponential inequality with two rates rather than one, and as such
would then give us a pretty complete understanding of this topic.

\end{document}